\documentclass[12pt]{amsart}
\usepackage{graphicx}
\usepackage{latexsym}
                 \def\s{\sigma}

       \def\sig{\hat{\s}}  
     \def\sn{\s^{-1}}
      
\def\Xi{\Xi}     \def\Pi{\Pi}       
             
    \def\L{\mathcal{L}}
\def\B{\mathcal{B}}    \def\O{\mathcal{O}}
\def\C{\mathcal{C}}    
    \def\S{\mathcal{S}} 
    \def\T{\mathcal{T}}
\newtheorem{Thm}{Theorem}[section]
\newtheorem{Cor}[Thm]{Corollary}
\newtheorem{Conj}[Thm]{Conjecture}
\newtheorem{Prop}[Thm]{Proposition}
\newtheorem{Lem}[Thm]{Lemma}

\newtheorem{Def}[Thm]{Definition}
\newtheorem{condition}[Thm]{Condition}
\newtheorem{formulla}[Thm]{}
\newtheorem{example}[Thm]{Example}
\newtheorem{problem}[Thm]{Problem}

\title{The fourth Skein module and the Montesinos-Nakanishi
       Conjecture for $3$-algebraic links}
\address{Dept. of Mathematics, the George Washington University, 
         Washington, DC 20052, USA \\ and University of Maryland,
         College Park}  \email{przytyck@gwu.edu}
\address{Dept. of 
        Mathematics, the George Washington University, Washington, 
        DC 20052, USA} \curraddr{Dept. of Mathematics, SUNY at 
        Buffalo, NY 14214, USA} \email{tsuka@tkg.att.ne.jp} 
\begin{document}\baselineskip=15pt\maketitle\pagestyle{plain}
%
\vspace{-3mm}\begin{center}
J\'{o}zef H. Przytycki\footnote[1]{partially supported by NSF grant 
DMS-9808955} and Tatsuya Tsukamoto\footnote[2]{partially supported 
by JSPS Research Fellowships for Young Scientists} \par\vspace{3mm}
In memory of Rodica Simion\footnote[3]{She was a member of the defence 
committee for the doctoral degree of the second author and she died 
just before the defence.} (January 18, 1955 - January 7, 2000)
\end{center}
%
%
\begin{abstract} We study the concept of the fourth skein module of 
$3$-manifolds, that is a skein module based on the skein relation 
$b_0L_0 + b_1L_1 + b_2L_2 +  b_3L_3 = 0$ and a framing relation 
$L^{(1)} = a L$ ($a$, $b_0$, $b_3$ invertible). We give necessary 
conditions for trivial links to be linearly independent in the module. 
We investigate the behavior of elements of the skein module under
the $n$-move and compute the values for $(2,n)$-torus links and twist 
knots as elements of the skein module. Using the idea of mutants and 
rotors, we show that there are different links representing the same 
element in the skein module. We also show that algebraic links 
(in the sense of Conway) and closed $3$-braids are linear combinations 
of trivial links. We introduce the concept of $n$-algebraic tangles 
(and links) and analyze the skein module for $3$-algebraic links. As 
a byproduct we prove the Montesinos-Nakanishi $3$-moves conjecture 
for $3$-algebraic links (including $3$-bridge links). In the case of 
classical links (i.e. links in $S^3$) our skein module suggests three 
polynomial invariants of unoriented framed (or unframed) links. One of 
them generalizes the Kauffman polynomial of links and another one can 
be used to analyze amphicheirality of links (and may work better than 
the Kauffman polynomial). In the end, we speculate about the meaning 
and importance of our new knot invariants. \end{abstract}
%
%
\section{Introduction}
To understand the structure of $3$-manifolds and links inside them, 
the first author introduced the concept of the skein module \cite{Pr-2} 
motivated by the Conway idea of ``linear skein". Skein modules are 
quotients of a free module of formal linear combinations of links in a 
$3$-manifold by properly chosen (local) skein relations. In search of 
appropriate skein relations, we analyze deformations of moves on links 
previously studied in knot theory. The simplest of such moves are $n$ 
(twist) moves (Figure $1$). The simplest skein modules are 
deformations of $1$-moves, \cite{Pr-4} (a linear skein relation). The 
Conway type skein modules are deformations of $2$-moves (a quadratic 
skein relation) and they are studied extensively (e.g. 
\cite{Pr-2,Tu,H-P,Bu,P-S}). In this paper we start a systematic study of 
skein modules based on a deformation of $3$-moves (a cubic skein 
relation). Since there are four terms involved in the relation, we call 
such a skein module the fourth skein module and denote it by $\S_4(M)$, 
following the notation of \cite{Pr-2}. 
In Section \ref{SecMNC}, we discuss the Montesinos-Nakanishi conjecture 
on $3$-moves. This conjecture, and its partial solutions, make the study 
of fourth skein modules feasible. 
In Section \ref{SecFSM}, we define the fourth skein module and we make 
the general conjecture about generators of the fourth skein module in 
$S^3$ and $B^3$ with $n$ boundary points ($n$-tangles).
In Section \ref{SecLITL}, we discuss the necessary conditions for 
coefficient ring so that trivial links are linearly independent in 
the fourth skein module. We prove the independence in two, rather 
degenerated cases. 
In Section \ref{SecGFSM}, we consider the case in which trivial links 
differ only by a multiplicative constant. We show that after proper 
substitution the skein module reduces to the Kauffman skein module, 
at least for links generated by trivial links. Generally we have 
always epimorphism into the Kauffman skein module.
In Section \ref{SecTorusandTwists}, we give an example of computation, 
finding eigenvectors of Dehn twists on $2$-tangles and applying this 
to describe the effect of $n$-moves on elements of the skein module. 
As a corollary we compute the values for $(2,n)$-torus links and twist 
knots.
In Section \ref{SecMRFSM}, we use the idea of mutants and rotors to 
construct different links representing the same element in the 
fourth skein module.
In Section \ref{SecnAlgTL}, we introduce a notion of an $n$-algebraic 
tangle and link, generalizing the notion of an algebraic link in the 
sense of Conway. We show that $n$-bridge links and closed $n$-braids 
are $n$-algebraic. We prove Conjectures \ref{GeneralConjMN} and 
\ref{ConjGenerator} for $2$- and $3$-algebraic tangles. 
In the last section, we speculate on existence of ``fourth" polynomial 
link invariants of classical links (in $S^3$) and show connection of 
our work with investigation of Coxeter, Assion and Wajnryb on finite 
quotients of braid groups.
%
%
%
\section{The Montesinos-Nakanishi conjecture}\label{SecMNC}
%
\hspace{2mm}
Let $L$ be a link (possibly a relative link) in a $3$-manifold $M$ 
(considered up to ambient isotopies). An $n$-move is a local change 
of a link which adds $n$ positive or negative half-twists to $L$ 
($L_0 {\longleftrightarrow}^{n-move} L_n$). In an oriented 
$3$-manifold we may distinguish  a move $L_0 \to L_n$ from a move 
$L_n \to L_0$. We occasionally use the notation $+n$-move and
$-n$-move in these cases (see Figure $1$).  
The first part of Conjecture \ref{GeneralConjMN} 
(Montesinos-Nakanishi conjecture) motivated, in part, this work. 
%
\begin{Conj}\begin{enumerate}\label{GeneralConjMN}
\item[$($i$)$]   $[Montesinos-Nakanishi$ $conjecture;$ \cite{Ki}$]$ 
                 Every link is $3$-equivalent to a trivial link 
                 $($i.e. it can be converted to a trivial link 
                 by $3$-moves$)$.
\item[$($ii$)$]  $($\cite{Ki}$)$ Every $2$-tangle is $3$-equivalent 
                 to a tangle with no more than one crossing 
                 (equivalently to one of the four ``basic" 
                 $2$-tangles, with possible trivial components); 
                 see Figure $2$.
\item[$($iii$)$] Every $3$-tangle is $3$-equivalent to one of the 
                 $40$ ``basic" $3$-tangles, with possible trivial 
                 components; see Figure $3$. Note that every basic 
                 $3$-tangle has no more than four crossings.
\item[$($iv$)$]  There is a finite number, $g(n)$, of ``basic" 
                 $n$-tangles, such that every $n$-tangle is 
                 $3$-equivalent to one of the $g(n)$ ``basic" 
                 $n$-tangles, with possible trivial components.
\item[$($v$)$]   $g(n)$ $=$ $\prod_{i=1}^{n-1}(3^i+1)$ $($e.g. 
                 $g(4)=1120)$.
\end{enumerate}\end{Conj}
%
\par\vspace{3mm}\begin{center} 
\framebox[2cm]{Figure $1$}\par\vspace{2mm}
\framebox[2cm]{Figure $2$}\par\vspace{2mm}
\framebox[2cm]{Figure $3$}\par\vspace{3mm}\end{center}
%
\section{The fourth skein module, $\S_4(M)$}\label{SecFSM}
%
In this section we start a systematic study of the fourth skein
module based on a linear relation among links $L_0,L_1, L_2$ and
$L_3$ of Fig.4 (the module was only shortly mentioned in \cite{Pr-6}).
\begin{Def}\label{DefFSM}
Let $M$ be an oriented $3$-manifold, $R$ a commutative ring with 
identity, and $\L_{fr}$ the set of ambient isotopy classes 
of unoriented framed links in $M$. Let $R\L_{fr}$ be the free 
$R$ module generated by $\L_{fr}$ and $a,b_0,b_1,b_2,b_3,$ 
elements in $R$ such that $a,b_0$ and $b_3$ are invertible. We 
define the {\rm fourth skein module}, $\S_4(M;R,a,b_0,b_1,b_2,b_3)$, 
or shortly $\S_4(M)$ as the quotient of $\L_{fr}$ by a 
submodule generated by the framing relations $L^{(1)}-aL$ and the 
fourth skein expression: $b_0L_0 + b_1L_1 + b_2L_2 + b_3L_3$.
Similarly, we can define a {\rm relative fourth skein module} 
starting from an oriented manifold $M$ with $2n$ points $x_1$, $x_2$,
$\ldots$, $x_{2n}$ chosen on $\partial M$ and the set of ambient 
isotopy classes of unoriented relative framed links in 
$(M, \partial M)$ such that $L \cap \partial M$ $=$ $\partial L$ 
$=$ $\{x_i\}$ for each relative link $L$. \end{Def}
%
\par\vspace{3mm}\begin{center}\framebox[2cm]{Figure $4$}
\par\vspace{3mm}\end{center}
%
Note that the following conjecture is a generalization of Conjecture 
\ref{GeneralConjMN}. In fact, Conjecture \ref{ConjGenerator} reduces 
to Conjecture \ref{GeneralConjMN} for $b_1=b_2=0$ and $b_0=-b_3$.
%
\begin{Conj}\begin{enumerate}\label{ConjGenerator}
\item[$($i$)$]  $\S_4(S^3)$ is generated by trivial links.
\item[$($ii$)$] The fourth skein module of $2$-tangles in a disk is 
                generated by tangles with no more than one crossing.
\item[$($iii$)$]The fourth skein module of $3$-tangles in a disk
                is generated by the $40$ basic $3$-tangles described 
                in Figure $3$ with possibly trivial components.
\item[$($iv$)$] There is a function $h(n)$ such that the fourth skein 
                module of $n$-tangles in a disk is generated by tangles
                with no more than $h(n)$ crossings.
\end{enumerate}\end{Conj}
\noindent
In this paper, we prove the conjecture for $3$-algebraic tangles and 
$3$-algebraic links (including $3$-bridge links and closures of 
balanced $3$-tangles). The notion of an $n$-algebraic tangle
(and link) is a new concept defined for the first time in this paper,
Section 8.
We end this section by listing several useful 
properties of the fourth skein module, including the {\it Universal 
Coefficient Property}, which our skein module shares with other skein 
modules \cite{Pr-2,Pr-6}.
%
\begin{Thm}\label{ThmUCP}\begin{enumerate}
\item [$(1)$] An orientation preserving embedding of $3$-manifolds
              $i: M \to N$ yields a homomorphism of skein modules 
              $i_{*}: \S_4(M) \to \S_4(N)$. The above
              correspondence leads to a functor from the category 
              of $3$-manifolds and orientation preserving embeddings 
              $($up to ambient isotopy$)$ to the category of 
              $R$-modules $($with specified elements $a,b_0,b_1$,
              $b_2,b_3\in R$, $a,b_0,b_3$ invertible$)$.
\item [$(2)$] $($Universal Coefficient Property$)$\\ 
              Let $r: R \to R'$ be a homomorphism of rings 
              $($commutative with $1)$. We can think of $R'$ as an 
              $R$ module. Then the identity map on $\L_{fr}$ 
              induces the isomorphism of $R'$ $($and $R)$ modules: 
              $$ \bar r: \S_4(M;R,a,b_0,b_1,b_2,b_3)
              \otimes_{R} R' \to \S_4(M;R',r(a),r(b_0),
              r(b_1), r(b_2),r(b_3)) .$$ In particular 
              $\S_4(M;Z,1,0,0,-1)= \S_4(M;Z[x],1,x,-x,-1) 
              \otimes_{Z[x]}Z$, where $r(x)=0$.
\item [$(3)$] Let $M=F\times I$ where $F$ is an oriented surface.
              Then $\S_4(M)$ is an algebra, where $L_{1} \cdot
              L_{2}$ is obtained by placing $L_{1}$ above $L_{2}$ 
              with respect to the product structure. The empty link 
              $T_{0}$ is the neutral element of the multiplication. 
              Every embedding $i:F' \to F$ yields an algebra 
              homomorphism $i_{*}: \S_4(F'\times I) \to
              \S_4(F\times I)$.
$\Box$\end{enumerate}\end{Thm}
%
%
%
\section{Linear independence of trivial links}\label{SecLITL}
%
\vspace{3mm}\hspace{2mm}
We analyze the general question: for which $a$ and $b_i$ are trivial 
links in $S^3$ linearly independent in $\S_4(S^3)$? We give necessary 
conditions and conjecture that they are also sufficient. From the 
relation $b_0L_0 + b_1L_1 +b_2L_2 + b_3L_3 =0$ (Figure $5$),  we get the 
relation $(b_0 + a^{-1}b_1 + a^{-2}b_2 + a^{-3}b_3)L_0 =0$ for any 
framed link $L_0$ (including the trivial knot). Thus to have the 
linear independence of trivial links we should assume:
%
\begin{condition}\label{CondtnTwistRelation}
$b_0 + a^{-1}b_1 + a^{-2}b_2 + a^{-3}b_3=0$
\end{condition} 
%
\par\vspace{3mm}\begin{center}\framebox[2cm]{Figure $5$}
\par\vspace{3mm}\end{center}
%
\noindent Now consider two ambient isotopic diagrams $P_1$ and $P_2$ 
(Figure $6$). Let $T_n$ be the trivial $n$ component link. For $P_1$ 
we compute: $P_1 = -b_3^{-1}(ab_2 + a^{-1}b_0 +b_1T_1)P$ 
and for $P_2$ we have: $P_2 = -b_0^{-1}(a^{-1}b_1 + ab_3 +b_2T_1)P$. 
Since $P_1 = P_2$ and $P$ can be any link, therefore in order 
to have linear independence we require: 
%
\begin{condition}\label{CondtnChainRelation}
\begin{enumerate}
\item[$(a)$] $b_0b_1=b_2b_3$, and
\item[$(b)$] $ab_0b_2 + a^{-1}b_0^2 = a^{-1}b_1b_3 + ab_3^2$
\end{enumerate}
\end{condition}
%
\par\vspace{3mm}\begin{center}\framebox[2cm]{Figure $6$}
\par\vspace{3mm}\end{center}
%
Condition \ref{CondtnTwistRelation} can be rewritten as (using 
Condition \ref{CondtnChainRelation} $(a)$ to eliminate $b_1$):
$$ a^3 + \frac{b_3}{b_0} + \frac{b_2}{b_0}(a^2\frac{b_3}{b_0} + a) 
= 0$$ and Condition \ref{CondtnChainRelation} $(b)$ can be written 
as (again using Condition \ref{CondtnChainRelation} $(a)$
to eliminate $b_1$): $$1 - a^2 (\frac{b_3}{b_0})^2 + 
\frac{b_2}{b_0}(a^2 - (\frac{b_3}{b_0})^2) = 0$$
If we eliminate $b_2$, then we get:
%
\begin{condition}\label{CondtnZeroDivisor}
$(a^4-1)(b_3^3 + ab_0^3)=0$
\end{condition}
\noindent
It is convenient to work with a ring without complicated zero divisors 
so we consider two cases of Condition \ref{CondtnZeroDivisor} separately.
%
\begin{enumerate}
\item[$(1)$]    Assume $a^4=1$. Our conditions reduce now to
                $(b_3+a^{3}b_0)(b_2+a^2b_0)=0$ which again 
                leads to two cases:
\begin{enumerate}
\item[$($i$)$]  Let $b_3 = - a^{3}b_0$. Then $b_2$ is a free variable, 
                and $b_1=-a^{-1}b_2$. As $b_0$ is invertible and 
                relations are homogeneous, we can put $b_0=1$.
                Let also write $b_1=ax$; then $b_2=-a^2x$, 
                $b_3=-a^{3}$. We work in this case with a ring 
                $R=Z[x][a]/(a^4-1)$. In particular we measure 
                framing only mod 4. Our skein relation has a form:
                $L_0 + axL_1 - a^2xL_2 -a^{3}L_3 = 0$.
\item[$($ii$)$] Let $b_2= -a^2b_0$. Then $b_3=-a^2b_1$ and after 
                substituting $b_0=1, b_1=ax$ one obtains the skein 
                relation: $L_0 + axL_1 - a^2L_2 -a^{3}xL_3 = 0$.
\end{enumerate}
\item[$(2)$]    Assume $b_3^3 + ab_0^3 = 0$. This leads to the 
                following solution $($for simplicity we put $b_0=1$ 
                and $b_3=b$$)$: $b_1=b(b^2+b^{-2})$, $b_2 = b^2+b^{-2}$, 
                and $a=-b^3$. The skein relation has now the form 
                $L_0 + b(b^2+b^{-2})L_1 + (b^2+b^{-2})L_2 +bL_3$, 
                and the framing relation $L^{(1)}=-b^3L$. $R=Z[b]$.
\end{enumerate}
%
We have considered conditions which are necessary for trivial links 
to be linearly independent in the fourth skein module, $\S_4(M)$. 
We conjecture that these conditions are also sufficient.
In the case of Condition \ref{CondtnZeroDivisor} $($1$)$$($ii$)$ 
we are able to prove it, and the skein module of $S^3$ seems 
to be rather trivial.
%
\begin{Thm}\label{ThmLinIndep}\begin{enumerate}
\item[$($i$)$]   If $a^4=1$, $b_2 = -a^2b_0$, $b_3 = -a^2b_1$ then 
               trivial links, $T_i$, are linearly independent in 
               $\S_4(S^3)$.
\item[$($ii$)$]  Consider a homomorphism $h: R\L_{fr} \to$ 
               $Z[x,t,a]/(a^4-1)$ given by $h(L) = a^{fr(L)}$ 
               $t^{com(L)}$ where $com(L)$ is the number
               of components of a link $L$ and $fr(L)$ is for 
               oriented framed link a difference between framing 
               of $L$ and $0$-framing in $S^3$. For unoriented 
               links $fr(L)$ is well defined $mod(4)$.
               Then $h$ yields a homomorphism $\hat h$ from the 
               fourth skein module described by Condition 
               $\ref{CondtnZeroDivisor}$ $(1)$$(ii)$  
               $($i.e. $b_0=1, b_1=ax, b_2= -a^2, b_3=-a^3x)$   
               to $Z[x,t,a]/(a^4-1)$. If we allow the empty link, 
               $T_0$, in $\L_{fr}$ and Conjecture 
               $\ref{ConjGenerator}$ $($i$)$ holds then $\hat h$ 
               is an algebra isomorphism, where product of two 
               links is defined to be their disjoint sum.
\end{enumerate}\end{Thm}
%
\begin{proof} To show $($ii$)$ it suffices to show that $h$ sends 
the skein relation $L_0 + axL_1 - a^2L_2 -a^{3}xL_3$ to $0$. 
We have $h(L_0 + axL_1 - a^2L_2 -a^{3}xL_3) = h(L_0  - a^2L_2) +
axh(L_1 - a^{2}L_3) = 0$ because $com(L_i)=com(L_{i+2})$ and
$fr(L_i)\equiv fr(L_{i+2}) + 2\ mod(4)$. The last statement 
follows from the fact that if $L$ is an oriented framed link 
and $D_L$ its diagram such that the framing of $L$ is the flat 
framing of $D_L$ then $fr(L) = Tait(D_L)= \Sigma_p sgn(p)$ where 
the sum is taken over all crossings of $D_L$. $h$ and $\hat h$ 
are clearly algebra epimorphisms and $h^{-1}$ and $\hat h^{-1}$, 
where $h^{-1}(t) = T_1$ are left inverses of $h$ and $\hat h$ 
respectively. Thus $\hat h$ is an isomorphism on the subspace of 
$\S_4(S^4)$ generated by trivial links. In fact our proof works 
for any rational homology sphere, and assuming  $a^2=1$, for any 
$3$-manifold. $($i$)$ follows immediately from $($ii$)$, because $t^i$ 
are linearly independent in the ring of polynomials.
\end{proof} \par
%
\begin{Conj}\label{ConjLinIndep}
The trivial links, $T_i$, are linearly independent in $\S_4(S^3)$ 
in the following cases: 
\begin{enumerate}
\item[$($i$)$]   $a^4=1$, $b_3=-a^3b_0$ and $b_1=-a^3b_2$   
\item[$($ii$)$]  $ab_0^3 = -b_3^3$, $b_2b_0^{-1} = b_3^2b_0^{-2} + 
                 b_3^{-2}b_0^2$ , $b_0b_1= b_2b_3$
\end{enumerate}\end{Conj}
\par\noindent
Notice that in all cases of Conjecture \ref{ConjLinIndep} and in 
Theorem \ref{ThmLinIndep}, we assumed that $b_0b_1=b_2b_3$ (compare 
Section \ref{SecGFSM}). Conjecture \ref{ConjLinIndep} $($1$)$$($i$)$ 
leads to the polynomial invariant of unframed links ($a=1$) in $S^3$ 
(for links generated by trivial links), $\S_4(L)(x,t)$. 
If $L= \Sigma_iw_i(x)T_i$ in our skein module, then 
$\S_4(L)(x,t)= \Sigma_iw_i(x)t^i$. Conjecture \ref{ConjLinIndep} 
$($ii$)$ leads to the polynomial invariant of framed links in 
$S^3$ (for links generated by trivial links), $\S_4(L)(b,t)$. 
If $L= \Sigma_iv_i(b)T_i$ in our skein module, then 
$\S_4(L)(b,t)= \Sigma_iv_i(b)t^i$. For $n$-tangles the analogue of 
Conjecture \ref{ConjLinIndep} is more involved but very interesting. 
We will consider below the case $($i$)$. \par
%
\begin{Conj}\label{ConjRelSM}\begin{enumerate}
\item[$(1)$] Consider the skein module $S_4(D(n),x)$ with $a=1$ 
             and the skein relation $L_0 +xL_1 -xL_2 -L_3$ where 
             $D(n)$ is a disk with $2n$ boundary points. Let 
             $B(n)$ be a set of $n$-tangles, one from each 
             $3$-move equivalence class. Then $B(n)$ is a base 
             for the fourth skein module $S_4(D(n),x)$, for
             $x\neq 1$.
\item[$(2)$] Consider the $Z_3$-linear space of all $3$-colorings 
             of a tangle, $b$, denoted by $Tri(b)$ and the 
             homomorphism $\phi_b: Tri(b)  \to Z_3^{2n}$ where 
             $ Z_3^{2n}$ is the space of all colorings of 
             boundary points \cite{Pr-5}. Then different elements 
             of $B(n)$ can be distinguished by their $3$-colorings. 
             Precisely, if $b_1\neq b_2$ then either 
             $rank (Tri(b_1)) \neq rank (Tri(b_2))$ or ranks are 
             the same but $\phi_{b_1}(Tri(b_1))
             \neq \phi_{b_2}(Tri(b_2))$.
\end{enumerate}\end{Conj} \par\noindent
%
By definition, $(1)$ holds for $x=0$, and so the theorem says here 
that the deformation ($0 \to x$) does not change the module. We 
exclude the case $x=1$ in the conjecture because we know 
that our skein module of $S^3$ (or of $D(0)$) behave differently 
for $x=0$ and $x=1$; in the first case we identify $3$-move 
equivalence classes and in the second case (as far as we know, 
compare Theorem \ref{ThmLinIndep}) links with the same number 
of components. It is known that each $3$-move is preserving a 
$3$-coloring of a tangle and we conjecture the inverse holds, 
that is, $n$-tangles with the same $3$-coloring structure are 
related by $3$-moves. Conjecture 4.6 suggests the importance of 
analyzing $Tri(b)$ and $\phi_{b}(Tri(b))$. We have proved recently 
that $\phi_{b}(Tri(b))$ are Lagrangians in the space of colorings 
of the boundary of the tangle $b$ with respect to the properly 
chosen symplectic structure \cite{Pr-7}.
%
\section{Generic case of the fourth skein module;
$b_0b_1 \neq b_2b_3$.}\label{SecGFSM}
\hspace{2mm}
In the previous section we showed that in all cases in which $T_i$ 
are (conjecturally) linearly independent, one has $b_0b_1 = b_2b_3$. 
We show that if $b_0b_1 - b_2b_3$ is invertible then one has:
\par\vspace{5mm}
%
\begin{formulla}\label{FormulaLinkwithTrivial}
$$L\sqcup T_1 = 
\frac{a^{-1}b_1b_3 -ab_0b_2 + ab_3^2 -a^{-1}b_0^2 }{b_0b_1 - b_2b_3}L$$
\end{formulla}\par
%
\noindent In particular: \par
%
\begin{formulla}\label{FormulaTrivialLinks}
$$T_n =(\frac{a^{-1}b_1b_3 -ab_0b_2 + ab_3^2 -a^{-1}b_0^2 }
{b_0b_1 - b_2b_3})^{n-1}T_1.$$
\end{formulla}\par\noindent
%
Namely the relation given by reversing the Hopf link summand gives:
$-b_3^{-1}(ab_2L + b_1(L\sqcup T_1) + a^{-1}b_0L)=
-b_0^{-1}(a^{-1}b_1L + b_2(L\sqcup T_1) + ab_3L)$. Thus
$(b_0b_1 - b_2b_3)(L\sqcup T_1) =$ $(a^{-1}b_1b_3 -ab_0b_2 + 
ab_3^2 -a^{-1}b_0^2)L$. This gives the formulas 
\ref{FormulaLinkwithTrivial} and \ref{FormulaTrivialLinks}.
We should also remember that we have found the restriction 
(Condition \ref{CondtnTwistRelation}):
$b_0 + a^{-1}b_1 + a^{-2}b_2 + a^{-3}b_3$ $=$ $0$. For $a^2\neq 1$ 
and the substitution $b_0=1,b_1=-(z+a),b_2=za+1,b_3=-a$
our skein module (polynomial) reduces to the Kauffman 
polynomial for links generated by trivial links (we get in 
particular $T_2=\frac{a+ a^{-1} -z}{z}T_1$). More precisely, 
our substitutions give the skein relations:\par
%
\begin{formulla}\label{Formula4SRrelKSM}
$L_0 -(z+a)L_1 +(za + 1)L_2 -aL_3 = 0$ , $L^{(1)}=aL$.
\end{formulla}\par\noindent
%
On the other hand the Kauffman relations are:\par
%
\begin{formulla}\label{FormulaKSRelation} 
$L_- + L_+ = z(L_0 + L_{\infty})$, $L^{(1)}=aL.$ 
\end{formulla}\par\noindent
%
We work then with the coefficients in a commutative ring with 
identity, $R$, where $z$ and $a$ are chosen invertible elements in 
$R$. The Kauffman skein module $\S_{3,\infty}(M;R,z,a)$ is defined 
to be $R\L_{fr}/$({\it Kauffman relations}). As a consequence of 
Kauffman relations we have $L\sqcup T_1 = \frac{a+ a^{-1} -z}{z}L$ 
in $\S_{3,\infty}(M;R,z,a)$. By applying \ref{FormulaKSRelation} 
twice we get: 
$$L_0 + L_2 = z(L_1 + a^{-1}L_{\infty}),$$
$$L_1 + L_3 = z(L_2 + a^{-2}L_{\infty}).$$
From this we get the relation: 
$$L_0 + L_2 - a(L_1 + L_3) = z(L_1 -aL_2),$$
which is equivalent to \ref{Formula4SRrelKSM}. Thus we have an 
epimorphism from the fourth skein module onto the Kauffman skein 
module (Proposition \ref{PropEpitoKSM}).\par
%
\begin{Prop}\label{PropEpitoKSM}
Let $M$ be an oriented $3$-manifold, $R$ a commutative ring 
with identity, and $a, z$ invertible elements in $R$. Then 
\begin{enumerate}
\item[$(a)$] We have an epimorphism 
             $$\phi: \S_4(M;relations\ \ref{Formula4SRrelKSM}) 
             \to \S_{3,\infty}(M;R,z,a)$$
\item[$(b)$] For $a^2-1$ invertible in $R$ and $M=S^3$, the 
             epimorphism $\phi$ restricted to the subspace generated 
             by trivial links is a monomorphism. If Conjecture 
             $\ref{ConjGenerator}$ $($i$)$ holds then $\phi$ is an 
             isomorphism.
\end{enumerate}\end{Prop}\par\noindent
%
\begin{proof}
Part $(a)$ follows from the fact that relations \ref{Formula4SRrelKSM} 
follow from the Kauffman relations, \ref{FormulaKSRelation}. If the 
trivial link $T_1$ is linearly independent in $\S_{3,\infty}(M;R,z,a)$ 
as is the case for $M=S^3$ or more generally for $M= F\times I$ (i.e. 
the product of a surface and the interval), and as is conjectured for 
any $3$-dimensional manifold, then we have a monomorphism $\psi$ from 
$RT_1$, the submodule generated by $T_1$ to $\S_4(M;relations$ 
\ref{Formula4SRrelKSM}$)$, which is left inverse to $\phi$ ($\psi \phi 
= Id_{RT_1}$). As observed before, for $a^2-1$ invertible $L\sqcup T_1 
= \frac{a+ a^{-1} -z}{z}L$ in the fourth skein module, thus $(b)$ of 
the proposition follows. \end{proof}  \par
%
\section {The fourth skein module for $(2,n)$-torus links and twist 
knots}\label{SecTorusandTwists}
\hspace{2mm}
We analyze in this section the behavior of the value of a link as
an element of the fourth skein module $\S_4(S^3)$ under $n$-moves 
and as a corollary compute the values of $(2,n)$-torus links and 
of twist knots. We show an example in the case of the skein relation 
$L_0 + axL_1 -a^2xL_2 -a^3L_3 =0$ (we discussed the importance of 
this relation, for $a^4=1$, in Section \ref{SecLITL}). 
\par\vspace{3mm}\noindent
%
Let $\tau$ denote the positive half twist on a tangle $L_0$,
taking framing into account we can write $\tau (L_0) =aL_1$, or
generally $\tau^k (L_i) =a^kL_{i+k}$. We will find the expression
for $a^nL_n$ in the basis $L_{-1}, L_0, L_1$. Of course there is
nothing sophisticated in our computation but it is useful to be 
able to evaluate our invariants for several classes of links.
Because $\tau (L_{-1}) = aL_0$, $\tau (L_0) = aL_1$, and 
$\tau (L_1) = aL_2 = aa^{-3}(-a^2xL_1 + axL_0 + L_{-1}) = 
a^{-2}L_{-1} + a^{-1}xL_{0} - xL_1 $, therefore the matrix, 
$A_{\tau}$ of our linear map has the form: \par\vspace{3mm}
%
$$A_{\tau} = \left (\begin{array}{ccc} 0 & 0 & a^{-2} \\
a & 0 & a^{-1}x \\ 0 &a & -x \end{array}\right ) $$ 
\par\vspace{3mm}\noindent
%
The characteristic polynomial of $\tau$ is therefore $\chi (A_{\tau})= 
-\lambda^3 - x\lambda^2 + x\lambda  +1 = (1- \lambda)(\lambda^2 + 
(x +1)\lambda + 1)$. As $\Delta = (x +1)^2 -4$, it is convenient to 
put $x+1 = -s^2 - s^{-2}$. Thus 1 is an eigenvalue of $\tau$ and two 
other eigenvalues are $\lambda_{\pm} = \frac{s^2 + s^{-2}\pm (s^2 - 
s^{-2})}{2}$, which gives $\lambda_+ = s^2$ and $\lambda_- = s^{-2}$. 
From this we get the eigenvectors of $\tau$ to be:     \par\vspace{1mm}
$E(1)=  a^{-1}L_{-1}  - (s^2 + s^{-2})L_0 + aL_1 $     \par\vspace{1mm}
$E(s^2)= a^{-1}s^{-1}L_{-1}  - (s + s^{-1})L_0 + asL_1$\par\vspace{1mm}
$E(s^{-2})= a^{-1}sL_{-1}  - (s + s^{-1})L_0 + as^{-1}L_1$
\par\vspace{1mm}\noindent   From these we obtain $L_0 =\frac{E(s^{-2}) 
-(s+s^{-1})E(1) + E(s^2)}{(s+s^{-1})(s-s^{-1})^2}$, and in general, 
we have the formula for $a^nL_n$. \par\vspace{1mm}
%
\begin{formulla}\label{FormulaTwisting}
$$a^nL_n =\tau^{n}(L_0) = \frac{s^{-2n}E(s^{-2}) -(s+s^{-1})E(1) + 
s^{2n}E(s^2)}{(s+s^{-1})(s-s^{-1})^2}$$ \par\vspace{1mm}\noindent
%
$$= \frac{s^{-2n}(a^{-1}sL_{-1} - (s + s^{-1})L_0 + as^{-1}L_1)} 
{(s+s^{-1})(s-s^{-1})^2}
  - \frac{(s+s^{-1}) (a^{-1}L_{-1} - (s^2 + s^{-2})L_0 + aL_1)}
{(s+s^{-1})(s-s^{-1})^2}$$ \par\vspace{1mm}\noindent
$$+ \frac{s^{2n}(a^{-1}s^{-1}L_{-1} - (s + s^{-1})L_0 + asL_1)}
{(s+s^{-1})(s-s^{-1})^2}
%
= \frac{(s^{2n-1}-(s+s^{-1}) + s^{-2n+1})a^{-1}L_{-1}}
{(s+s^{-1})(s-s^{-1})^2}$$ \par\vspace{1mm}\noindent
$$- \frac{(s+s^{-1})(s^{2n}-(s^2+s^{-2}) + s^{-2n})L_0}
{(s+s^{-1})(s-s^{-1})^2}
  + \frac{(s^{2n+1}-(s+s^{-1}) + s^{-2n-1})aL_{+1}}
{(s+s^{-1})(s-s^{-1})^2}$$ \par\vspace{1mm}\noindent
%
$$= \frac{(s^n-s^{-n})(s^{n-1}-s^{1-n})a^{-1}L_{-1}}
{(s+s^{-1})(s-s^{-1})^2}
  - \frac{(s+s^{-1})(s^{n-1}-s^{1-n})(s^{n+1} - s^{-n-1})L_0}
{(s+s^{-1})(s-s^{-1})^2}$$ \par\vspace{1mm}\noindent
$$+ \frac{(s^{n+1} - s^{-n-1})(s^n-s^{-n})aL_{+1}}
{(s+s^{-1})(s-s^{-1})^2}.$$ \end{formulla}\par\noindent
%
\begin{Cor}\label{CorDehnTwist} Consider the Dehn twist $\tau^2$ on 
the tangle. It has eigenvalues $1$, $s^4$, and $s^{-4}$ and the same 
eigenvectors as $\tau$. \end{Cor}
\noindent
One can use this observation to construct a torsion in a manifold
with an incompressible $S^2$ (e.g. $S^1 \times S^2$). In the spirit 
of \cite{Pr-1} we can check the condition at which an $n$-move 
changes our invariant only by a constant.
%
\begin{Cor}\label{CorTwisting}
If $s^{2n}=1$ and $s^4 \neq 1$ for $n$ even and  $s^2 \neq 1$ 
for $n$ odd then $ a^nL_n =\tau^{n}(L_0) = L_0$. In particular, 
as we know, if $x=-s^2-s^{-2}-1= -\frac{s^{3}-s^{-3}}
{s-s^{-1}} = 0$, $\tau^3(L_0)=a^3L_3 = L_0$. \end{Cor}
\noindent
To find the expression for torus links, $T_{2,n}$ (Figure $7$ ($a$)), 
in the fourth skein module, we observe that $aT_{2,-1} = a^{-1}T_{2,1}= 
T_1$ and $T_{2,0} =T_2$.
%
\par\vspace{3mm}\begin{center}\framebox[2cm]{Figure $7$}
\par\vspace{3mm}\end{center}
%
\begin{Cor}\label{CorTorus}
The invariant of the torus link of type $(2,n)$ is given by:
\begin{displaymath} a^nT_{2,n}=\tau^n(T_{2,0})
= \frac{(s^n-s^{-n})(a^{-2}(s^{n-1}-s^{-n+1}) + a^2(s^{n+1} - 
s^{-n-1}))}{(s+s^{-1})(s-s^{-1})^2}T_1 \end{displaymath}
\begin{displaymath}
- \frac{(s+s^{-1})(s^{n-1}-s^{-n+1})(s^{n+1} - s^{-n-1})}
{(s+s^{-1})(s-s^{-1})^2}T_2.
\end{displaymath}\end{Cor}
%
Notice that for unoriented links $T_{2,2}=T_{2,-2}$ which gives
$(a^2-a^{-2})((s^2+ 1 +s^{-2})T_2 + (-a^2 -a^{-2})T_1)=0$. We have 
already discussed the condition $a^4=1$ which is necessary if we
want trivial links to be linearly independent. The condition 
$(s^2+ 1 +s^{-2})T_2 = (a^2 +a^{-2})T_1 $ leads to 
$T_2 = \frac{a^2 +a^{-2} }{s^2+ 1 +s^{-2}}T_1 = 
-\frac{a^2 +a^{-2}}{x}T_1$ (this agrees with
\ref{FormulaTrivialLinks}). \par\vspace{2mm}\noindent
%
A similar calculation for twist knots, $a^n T[2,n]$ (Figure $7$ ($b$)) 
uses the Formula \ref{FormulaTwisting} and the initial conditions 
$T[2,0]=a^2T_1$, $T[2,-1] = a^{-1}T_1$, and \par $T[2,1]= T_{2,3}= 
a^{-3}\frac{(s^3-s^{-3})(a^{-2}(s^{2}-s^{-2}) + a^2(s^{4} - s^{-4}))T_1
-(s+s^{-1})(s^{2}-s^{-2})(s^{4} - s^{-4})T_2}{(s+s^{-1})(s-s^{-1})^2}.$
\par\vspace{3mm}
%
More generally we can write a straightforward formula for the pretzel 
link $P[n_1,n_2,...n_k]$, Figure $8$ by first using Formula 
\ref{FormulaTwisting} for each column to get a linear 
combination generated by torus links and then using the formula
again for these torus links to get a linear combination generated
by trivial links. To write a formula for a $2$-bridge link (and 
tangle) or for an algebraic tangle, we should first decompose 
an algebraic tangle into $n$-twist tangles and then, from the 
table of multiplication of the four basic $2$-tangles, we can read 
off the result. Using the table and this formula, we can also
figure out the automorphism of the ``rotation". A similar formula 
can be built for $3$-algebraic tangles (with knowledge of the
$40\times 40$ multiplication table, and the ``rotation" formula;
see Section \ref{SecnAlgTL}). 
%
\par\vspace{3mm}\begin{center}\framebox[2cm]{Figure $8$}
\par\vspace{3mm}\end{center}
%
\section{Mutations and rotations in the fourth skein module}
\label{SecMRFSM}
\hspace{2mm}
Using the idea of mutants and rotors, we show that there are different 
links representing the same element in the fourth skein module. If we 
have a $2$-tangle, $P$, we consider three involutions $m_x$, $m_y$,
and $m_z$ on it which denote the rotation of the tangle along the 
axis denoted in the subscript, Figure $9$. If we have a link diagram 
with a $2$-tangle in it, a mutation of the diagram along the tangle 
is obtained by performing one of the described involutions on the 
tangle (see Figure $9$). \par\vspace{3mm}
%
\par\vspace{3mm}\begin{center}\framebox[2cm]{Figure $9$}
\par\vspace{2mm}\framebox[2cm]{Figure $10$}
\par\vspace{3mm}\end{center}
%
\begin{Prop}\label{PropMutants}
Consider a $2$-tangle part of a link diagram. If it is generated by 
tangles of Figure $2$ $($with possible trivial 
components$)$, then any mutation of the link preserves the element of 
the fourth skein module.\end{Prop}
%
\begin{proof} Any tangle of Figure $2$ is symmetric with respect to 
$x$, $y$, and $z$ axes. \end{proof} \par\vspace{2mm}\noindent
%
Consider a tangle $R$ in a regular $n$-gon $D$ in $S^2$ ($n \geq 3$).
We call $R$ an ($n$-)rotor if $R$ is invariant under the rigid rotation
$\rho: D \rightarrow D$ by angle $2\pi/n$ around the center of $D$.
Take a link $L$ and let $L$ have a tangle decomposition $R \cup S$,
i.e., $L$ has a projection which intersects $\partial D$ 
transversely in $2$ points on each edge of $D$. If $R$ is an $n$-rotor,
then we call $S$ a stator. Consider a line which passes through the
center of $D$ and either a corner or middle point of an edge. Let
$\mu: D \rightarrow D$ be a $\pi-$rotation through the third dimension
with the line as an axis. Then we call $L' = \mu(R) \cup S$ a rotant
of $L$. Here note that we do not assume that $R$ is invariant under
the operation $\mu$. If we consider a similar rotation $\mu'$ for
a disk $S^2 - D$, then we have the same result $L' = R \cup \mu'(S)$.
See \cite{A-P-R} for more details. \par\vspace{2mm}\noindent
%
Regard each tangle of Figure $3$ as a stator of a $3$-rotor. Then 
$34$ of the $40$ basic $3$-tangles are invariant under an operation 
$\mu'$ for an axis of a symmetry of the polygon. Exceptions are 
$\s_1s_2$, $\s_1U_1\sn_1$, and their mirror images as they have only 
plane symmetries. Thus we have the following proposition. 
\par\vspace*{3mm}
%
\begin{Prop}\label{PropRotants}
Assume that a link $L$ has a tangle decomposition $R \cup S$, where
$R$ is a $3$-rotor and $S$ is a $3$-stator. If $S$ is generated by
$34$ tangles in Figure $3$ as described above (with possible trivial 
components), then its rotant and $L$ are the same in the fourth skein 
module. $\Box$\end{Prop}
%
\section{$n$-algebraic tangles and links}\label{SecnAlgTL}
%
\hspace{2mm}
In this section, we study the fourth skein module of the $3$-sphere
and the relative skein module of the 3-dimensional disk. 
To understand the structure better, we also study fourth skein 
algebra for $n$-braid groups. Here note that each skein module in 
this section has its tangle elements (including links as 
$0$-tangles) in $S^3$ or in $D^3$. Therefore we can analyze their 
diagrams on $\bf{R}^2$ ($\subset$ $S^2$) or $D^2$ and count the 
number of crossings of these diagrams. We use $\s_i$ for a crossing 
of a braid, as usual and use $\sig_i$ in the case we do not specify 
whether the crossing is positive or negative. We use three kinds of 
equality symbols to denote various relations between 
a pair of n-tangle diagrams (including link diagrams) $a$ 
and $b$. We say that $a$ and $b$ are {\it skein 
module equivalent}, denoted by $a$ $=$ $b$, if they are equal 
in the skein module. We say that $a$ and $b$ are {\it strongly 
equivalent}, denoted by $a$ $\equiv$ $b$, if they are ambient 
isotopic. Finally, we say that $a$ and $b$ are {\it dot equivalent}, 
denoted by $a$ $\doteq$ $b$, if the difference of $a$ and $b$ is 
presented in the skein module as a linear combination of elements 
with fewer crossings 
than $a$ and $b$ have. If every element of the linear combination 
has no more than $n$ crossings, then we denote the linear combination 
by $\O(n)$. \par
%
\subsection{$n$-algebraic tangles and $n$-algebraic links}
\label{SubSecNTNAL}
\hspace{2mm}
Let $\T$ denote the category of unoriented tangles and $\T(n)$ a 
semigroup of $n$-tangles. We ``read" morphisms from left to right 
and compose tangles in the same manner (it corresponds to the 
notation $P_1P_2$ denoting composition of tangles $P_1$ and $P_2$). 
$\T(n)$ allows the natural action by the $D_{2n} \oplus Z_2$ group, 
where $D_{2n}$ is the action of the dihedral group (preserving 
orientation) of $4n$ elements, generated by the rotation $r$ along 
the $z$-axis by the angle $2\pi/2n$, and the rotation $\pi$ along 
the $y$-axis denoted by $m_y$ and often called a mutation. The 
factor $Z_2$ corresponds to the mirror image map with respect to 
the $x$, $y$ plane. The semigroup of $n$-tangles has two important 
subsemigroups: the braid group $B_n$, generated by crossings 
(tangles) $\sigma_i$, and a semigroup of balanced $n$-tangles, 
$W(n)$, generated by $\sigma_i$ and tangles $U_i$ creating maxima 
and minima (for a visual notation, $)(^i$ might be good. See Figure 
$3$ for an example). In fact $U_1$ would suffice, as $U_{i+1} = 
\sigma_i\sigma_{i+1}U_i  \sigma_{i+1}^{-1}\sigma_i^{-1}$.
Representations of $W(n)$ (their finite dimensional quotients) were
studied by B.Westbury. We are introducing also another semigroup, 
generalizing Conway's algebraic tangles and allowing successful
induction.
%
\begin{Def}\label{DefnAlgTangles}\begin{enumerate}
\item[$(1)$] We define an n-algebraic tangle inductively as follows:
\begin{enumerate}
\item[$($i$)$]  An $n$-tangle with no more than one crossing is
                $n$-algebraic.
\item[$($ii$)$] Inductive step: If $P_1$ and $P_2$ are two 
                $n$-algebraic tangles then their twisted sum, 
                $r^i(P_1)r^j(P_2)$ is also $n$-algebraic.
\end{enumerate}
\item[$(2)$]  In a more restrictive way, we define $(n,k)$-algebraic
              tangles if in step $($ii$)$ $P_2$ has no more than $k$ 
              crossings.
\item[$(3)$]  An $n$-algebraic link $($resp. $(n,k)$-algebraic 
              link$)$ is a link with a diagram obtained by 
              closing an $n$-algebraic $($resp. $(n,k)$-algebraic$)$ 
              tangle by pairwise disjoint arcs.
\end{enumerate}\end{Def}
\noindent
We denote by $A_n$ the monoid generated by $n$-algebraic tangles.
%
\begin{Lem}\label{LemN1Algebraic}\begin{enumerate}
\item{(a)} Tangles in $B_n$, $W_n$ and $n$-bridge tangles are 
           $(n,1)$-algebraic.
\item{(b)} $n$-bridge links and closed n-braids are $(n,1)$-algebraic 
           links.
\end{enumerate}\end{Lem}
%
\begin{proof} The fact that elements of $B_n$ and $W_n$ are 
$(n,1)$-algebraic follows from the definition. For an $n$-bridge 
tangle, we start from the $n$-tangle which is equal to 
$U_1U_3...U_{2k-1}$ for $n=2k$ or $n=2k+1$, and then we add crossings 
one by one (compare Figure $11$, the figure is drawn to stress the 
fact that we deal with a 3-bridge tangle. In the algebra of 3-tangles 
we would describe it as $r^{-1}(r(U_1\sigma_2)\sigma_2)$). \end{proof}
%
\par\vspace{3mm}\begin{center}\framebox[2cm]{Figure $11$}
\par\vspace{3mm}\end{center}
%
Our main result in this section is the following.
%
\begin{Thm}\label{ThmThreeLinks} Any $3$-algebraic link can be 
generated as a linear combination of trivial links in $\S_4(S^3;R)$. 
\end{Thm}             
\noindent
Our proof works for any coefficients thus as a corollary we have:
%
\begin{Thm}\label{ThmM-NConjecuture} The Montesinos-Nakanishi 
conjecture is true for $3$-algebraic links. $\Box$\end{Thm}             
%
\noindent{\bf Remark.} What concerns the Montesinos-Nakanishi 
conjecture, we have two candidates for a possible counterexample. One 
is the double of Borromean rings proposed by Y. Nakanishi \cite{Ki}; 
it has 24 crossings. The other is proposed by Q. Chen in \cite{Ch} and 
it has $20$ crossings. Chen's link is a reduction of the closure of 
the 5-braid $(\s_1\s_2\s_3\s_4)^{10}$ by $3$-moves. 
In \cite{Ts-1,Ts-2}, it is shown that the double of Borromean rings 
and Chen's link are $4$-algebraic (in fact, $(4,6)$-algebraic) up to 
$3$-moves, which tell us that it is hard to prove Theorems 
\ref{ThmThreeLinks} and \ref{ThmM-NConjecuture} for $4$-algebraic 
links (it is shown in \cite{Ts-2} that these theorems hold for 
$(4,5)$-algebraic links).
%
\par\vspace{3mm}\begin{center}\framebox[2cm]{Figure $12$}
\par\vspace{3mm}\end{center}
%
\subsection{The fourth skein module for the $3$-algebraic links}
\label{SubSecFSM3Alg}
\hspace{2mm}
Each of the propositions and theorems presented in this section holds
also for $n = 2$ (we take $\{Id$, $\s_1^{\pm}\}$ and $\{U_1\}$ 
instead of $\B_3$ and $\C_3$, respectively). Since proofs are much
simpler in that case, we omit them. First of all, let us define the
fourth skein algebras for $A_n$ and $B_n$.
%
\begin{Def}\label{Def4SAforBG}
Let $R$ be commutative ring with identity. Let $RB_n$ be the 
$R$-module with the basis $B_n$. Furthermore, let $I_4$ be the ideal of 
$RB_n$ generated by the skein relation $b_0\s_i^0$ $+$ 
$b_1\s_i$ $+$ $b_2\s_i^2$ $+$ $b_3\s_i^3$, where $b_0$ and 
$b_3$ are invertible in $R$. We define the fourth skein module 
$\S_4(B_n;R)$ for $B_n$ as the quotient $RB_n/I_4$. The product 
for elements of $B_n$ induces a bilinear map $\S_4(B_n;R)$ $\times$ 
$\S_4(B_n;R)$ $\rightarrow$ $\S_4(B_n;R)$ so that, with respect 
to this product, $\S_4(B_n;R)$ becomes an algebra. Thus we call 
$\S_4(B_n;R)$ the fourth skein algebra for $B_n$. We define the
forth skein algebra $\S_4(A_n;R)$ for $A_n$ as a subalgebra
of the algebra of $n$-tangles modulo the fourth skein relations and
the framing relations generated by $n$-algebraic tangles.
\end{Def}\par 
%
Note that in $\S_4(B_3;R)$, any braid which contains $\sig_i^2$ or 
any non-alternating expression of $(\sig_1\sig_2)^2$ or 
$(\sig_2\sig_1)^2$ can be generated as a linear combination of 
elements with fewer crossings than those we start from. From the 
following proposition, we also see that $(\sig_1\sig_2)^2\sig_1$ 
and $(\sig_2\sig_1)^2\sig_2$ can be generated as a linear combination 
of elements with fewer crossings than those we start from. We call 
a braid in $\S_4(B_3;R)$ {\it reducible} if the braid contains one 
of those configurations. Otherwise, we call the braid {\it irreducible}.
%
\begin{Prop}\label{Prop3BraidofLength4}
Four configurations $(\s_1\sn_2)^2$, $(\sn_1\s_2)^2$, $(\s_2\sn_1)^2$, 
and $(\sn_2\s_1)^2$ are dot equivalent to each other in $\S_4(B_3;R)$.
\end{Prop}
%
\begin{proof}
It is sufficient to show that $(\sn_1\s_2)^2$ $\doteq$ $(\s_1\sn_2)^2$, 
$(\sn_2\s_1)^2$ $\doteq$ $(\s_1\sn_2)^2$, and $(\s_2\sn_1)^2$ $\doteq$ 
$(\sn_2\s_1)^2$ using the fourth skein relations.
\par\vspace{3mm}
$(\sn_1\s_2)^2$ $\equiv$ $\s_1^{-2}\sn_2\s_1\s_2^2$
$=$ $(-b_0^{-1}b_3)$$(-b_0b_3^{-1})$ $\s_1\sn_2\s_1\sn_2$
$+$ $(-b_0^{-1}b_3)$$(-b_2b_3^{-1})$ $\s_1\sn_2\s_1\s_2$
\par\vspace{3mm}\hspace{1cm}
$+$ $(-b_0^{-1}b_1)$$(-b_0b_3^{-1})$ $\sn_1\sn_2\s_1\sn_2$
$+$ $(-b_0^{-1}b_1)$$(-b_2b_3^{-1})$ $\sn_1\sn_2\s_1\s_2$
$+$ $\O(3)$
\par\vspace{3mm}\hspace{1cm}
$=$ $(-b_0^{-1}b_3)$$(-b_0b_3^{-1})$ $\s_1\sn_2\s_1\sn_2$ $+$ $\O(3)$.
\par\vspace{3mm}
$(\sn_2\s_1)^2$ $\equiv$ $\s_1\s_2^2\s_1^{-2}\sn_2$
$=$ $(-b_0^{-1}b_3)$$(-b_0b_3^{-1})$ $\s_1\sn_2\s_1\sn_2$
$+$ $(-b_0^{-1}b_3)$$(-b_2b_3^{-1})$ $\s_1\s_2\s_1\sn_2$
\par\vspace{3mm}\hspace{1cm}
$+$ $(-b_0^{-1}b_1)$$(-b_0b_3^{-1})$ $\s_1\sn_2\sn_1\sn_2$
$+$ $(-b_0^{-1}b_1)$$(-b_2b_3^{-1})$ $\s_1\s_2\sn_1\sn_2$
$+$ $\O(3)$
\par\vspace{3mm}\hspace{1cm}
$=$ $(-b_0^{-1}b_3)$$(-b_0b_3^{-1})$ $\s_1\sn_2\s_1\sn_2$ $+$ $\O(3)$.
\par\vspace{3mm}
$(\s_2\sn_1)^2$ $\equiv$ $\s_2^2\s_1\sn_2\s_1^{-2}$
$=$ $(-b_0^{-1}b_3)$$(-b_0b_3^{-1})$ $\sn_2\s_1\sn_2\s_1$
$+$ $(-b_0^{-1}b_3)$$(-b_2b_3^{-1})$ $\s_2\s_1\sn_2\s_1$
\par\vspace{3mm}\hspace{1cm}
$+$ $(-b_0^{-1}b_1)$$(-b_0b_3^{-1})$ $\sn_2\s_1\sn_2\sn_1$
$+$ $(-b_0^{-1}b_1)$$(-b_2b_3^{-1})$ $\s_2\s_1\sn_2\sn_1$
$+$ $\O(3)$
\par\vspace{3mm}\hspace{1cm}
$=$ $(-b_0^{-1}b_3)$$(-b_0b_3^{-1})$ $\sn_2\s_1\sn_2\s_1$ $+$ $\O(3)$.
\end{proof} \par\noindent
%
Define $\B_3$ as the set of $24$ invertible (braid type) basic 
tangles in Figure $3$. Then we have the following proposition.
%
\begin{Prop}\label{PropIrreducible3Braids}
For any irreducible element $b$ of $B_3$, there exists an element
$b'$ of $\B_3$ such that $b'$ is dot equivalent to $b$.
\end{Prop}                                                       
%
\begin{proof}
Let $n$ be the number of crossings of $b$. If $n$ $\leq$ $3$, then 
there exists an element $b'$ of $\B_3$ such that $b'$ is strongly 
equivalent to $b$. If $n$ $=$ $4$, then $b$ is an alternating 
expression of $(\sig_1\sig_2)^2$ or $(\sig_2\sig_1)^2$. Then, $b$ 
is dot equivalent to $(\s_1\sn_2)^2$ from Proposition 
\ref{Prop3BraidofLength4}. There is no irreducible element of $B_3$ 
with more than four crossings. \end{proof}\par\vspace{2mm}
%
\noindent Using this proposition, we obtain the following theorem.
%
\begin{Thm}\label{ThmThreeBraids} Any element of $\S_4(B_3;R)$ can be 
generated as a linear combination of elements of $\B_3$.\end{Thm}  
%
\begin{proof}
We prove the theorem by an induction on the number of crossings of 
the element. If $n = 0$, then clearly the statement holds. Assume 
that the statement holds in the case $n$ fewer than $k$. Consider 
the case $n = k$. If the element is reducible, then it is generated 
by elements with fewer crossings than $k$, which can be generated by 
elements of $\B_3$ from the assumption of the induction. If the 
element is irreducible, then it is generated by an element of $\B_3$ 
and elements with fewer crossings than $k$ from Proposition
\ref{PropIrreducible3Braids}. This induces that the irreducible element 
is also generated by elements of $\B_3$. \end{proof}\par\noindent
%
Define $\C_3$ as the set of $16$ non-invertible basic tangles in
Figure $3$. Then we have the following.
%
\begin{Prop}\label{PropProduct} A product of any pair of elements of 
$\B_3$ and $\C_3$ can be generated as a linear combination of elements 
of $\B_3 \cup \C_3$ (with possible trivial components). \end{Prop}
%
\begin{proof} If both elements of the pair are $3$-braids, then the 
statement follows from Theorem \ref{ThmThreeBraids}. The other cases
are easy to be checked. For simplicity, we ignore trivial components
in the consideration. \end{proof}
%
\begin{Prop}\label{PropRotation} Any rotation of elements of $\B_3$ 
and $\C_3$ can be generated as a linear combination of elements of 
$\B_3 \cup \C_3$. \end{Prop}
%
\begin{proof} We need to show only for $r^2((\s_1\sn_2)^2)$ and 
$r^{-1}((\s_1\sn_2)^2)$. From Figure $13$, we have \par\vspace{2mm}
$r^2((\s_1\sn_2)^2) = (\s_1\sn_2)^2 
              + b_0^{-1} b_2(\sn_2\sn_1 - \s_1 U_2\sn_1)
              + b_1   b_3^{-1}(\s_2\s_1 - \sn_1U_2\s_1)$
\par\vspace{3mm}\hspace{2cm}
             $+ b_0^{-1}b_1b_2b_3^{-1} (\sn_2\s_1  - \s_1 U_2\s_1)$
\par\vspace{3mm}
$r^{-1}((\s_1\sn_2)^2) = (\sn_1\s_2)^2
              + b_0^{-1}   b_2 (\sn_2\sn_1 - \s_1 U_2\sn_1)
              + b_1   b_3^{-1} (\s_2\s_1   - \sn_1U_2\s_1)$
\par\vspace{3mm}\hspace{2cm}
            $+ b_0^{-1}b_1b_2b_3^{-1} (\s_2\sn_1  - \sn_1U_2\sn_1)$
\end{proof}
%
\par\vspace{3mm}\begin{center}\framebox[2cm]{Figure $13$}
\par\vspace{3mm}\end{center}\par\vspace{2mm}\noindent
%
Since any $3$-tangle with no more than one crossing can be expressed
as $a^i t$, where $i = 0, \pm 1$, $t$ is an element of $\B_3
\cup \C_3$, and $a$ reflects a possible framing change yielded by
a ``kink", we obtain the following theorem from Propositions
\ref{PropProduct} and \ref{PropRotation}.
%
\begin{Thm}\label{ThmThreeTangles} Any element of $\S_4(A_3;R)$ can 
be generated as a linear combination of elements of $\B_3 \cup \C_3$
(with possible trivial components). $\Box$\end{Thm}     
%
\noindent {\it Proof of Theorem \ref{ThmThreeLinks}.} 
From Theorem \ref{ThmThreeTangles}, we know that every $3$-algebraic 
link can be generated as a linear combination of $3$-algebraic links 
obtained from elements of $\B_3 \cup \C_3$. For each of these links,
it is easy to see that it can be generated as a linear combination of 
trivial links. $\Box$
%
\section{Speculations}\label{SecSpeculations}
\hspace{2mm}
In this section we present results of calculations that suggest that 
even for $S^3$ the fourth skein module is more powerful (can 
distinguish more links) than  the third (Jones-Conway) and Kauffman 
skein modules. In particular, if the fourth skein polynomial exists 
(trivial links are linearly independent), then the polynomial is (at 
least sometimes) more powerful than the Jones-Conway (Homflypt) and 
Kauffman polynomials. Our calculation shows that under these 
assumptions we can  distinguish the $9_{42}$-knot from its mirror 
image $\bar 9_{42}$ (\cite{Ro}) and some 2-bridge links that share 
the same Jones-Conway and Kauffman polynomials (\cite{Ka-1,Ka-2}). 
In the second part of this section we present general ideas that can 
lead to a solution of Montesinos-Nakanishi conjecture and its 
generalizations. In particular we speculate about the geometrical 
meaning of the finite quotients of the braid group described by 
J.Assion \cite{A-1,B-W, Wa} (Symplectic and Unitary cases).
%
\begin{Conj}\label{}\begin{enumerate}
\item[$(1)$]   There is a polynomial invariant of unoriented links,
              $P_1(L) \in Z[x,t]$ that satisfies:
\begin{itemize}
\item[$($i$)$]   Initial conditions: $P_1(T_n) = t^n$, where $T_n$ is 
              the trivial link of $n$ components.
\item[$($ii$)$]  Skein relation, $P_1(L_0) + xP_1(L_1) - xP_1(L_2) 
              - P_1(L_3)=0$, where $L_0,L_1,L_2,L_3$ is a standard, 
              unoriented skein quadruple $(L_{i+1}$ is obtained 
              from $L_{i}$ by a right-handed half twist on two 
              arcs involved in $L_{i})$.
\end{itemize}
\item[$(2)$]   There is a polynomial invariant of unoriented framed 
              links, $P_2(L) \in Z[b^{\pm 1},t]$ which satisfies:
\begin{enumerate}
\item[$($i$)$]   Initial conditions: $P_2(T_n) = t^n$,
\item[$($ii$)$]  Framing relation: $P_2(L^{(1)}) = -b^3P_2(L)$, where 
              $L^{(1)}$ is obtained from a framed link $L$ by a 
              positive half twist on its framing.
\item[$($iii$)$] Skein relation: $P_2(L_0) + b(b^2 + b^{-2})P_2(L_1) +
              (b^2 + b^{-2})P_2(L_2) + bP_2(L_3)=0$.
\end{enumerate}
\item[$(3)$]   In each case the polynomial is uniquely defined. 
\end{enumerate}\end{Conj}
%
\begin{example}\label{10.1}
Consider the $9_{42}$-knot, in notation of \cite{Ro}. Then the 
polynomial $P_2(b,t)$ distinguishes $9_{42}$ from its mirror image 
$\bar 9_{42}$. In fact, we obtain the following for $9_{42}$.
We use the fact that we have the formula: $P_2(L)(b,t) 
= P_2(\bar L)(b^{-1},t)$. \end{example}
$P_2(9_{42})(b,t) =
(3b^{11}+7b^7+9b^3+8b^{-1}+6b^{-5}+4b^{-9}+3b^{-13}+b^{-17})$ $t$
$-$ $(b^{13}+6b^9+14b^5+20b$ \par\vspace{1mm}\hspace{10mm}
$+19b^{-3}+12b^{-7}+5b^{-11}+b^{-15})$ $t^2$
$+$ $(b^{11}+4b^7+8b^3+10b^{-1}+8b^{-5}+4b^{-9}+b^{-13})$ $t^3$.
%
\par\vspace{3mm}\begin{center}\framebox[2cm]{Figure $14$}
\par\vspace{3mm}\end{center}
%
\begin{example}\label{10.2}
Consider two $2$-bridge knots $K_1$ and $K_2$ as shown in the 
following figure. According to Kanenobu \cite{Ka-1,Ka-2,K-S}, they 
share the Jones-Conway $($Homflypt$)$ and Kauffman polynomials. We 
can, however, distinguish them by $P_1(x,t)$ polynomial. In fact, we 
obtain the following for $x = -2$; $P_1(K_1)(-2,t) = 49 t - 48 t^2$ 
and $P_1(K_2)(-2,t) = 28 t - 27 t^2$. \end{example}
%
\par\vspace{3mm}\begin{center}\framebox[2cm]{Figure $15$}
\par\vspace{3mm}\end{center}
%
We speculate that Conjectures \ref{GeneralConjMN} and 
\ref{ConjGenerator} can be approached by using results of Coxeter 
and of Assion and by interpreting them using some elementary moves 
on oriented links. Coxeter showed that $C_n= B_n/(\sigma_i)^3$ is 
finite iff $n \leq 5$, \cite{Co,A-2,Mu}. Assion found two basic 
cases in which $C_n/(Ideal)$ is finite: the ``symplectic" and 
``unitary" cases \cite{A-1,B-W,Mu}. Let $\Delta^5 = 
(\sigma_1\sigma_2\sigma_3\sigma_4)^5$ be a generator of the center 
of the braid group, $B_5$. It is easy to check that Assion's
ideals are generated by $\Delta^{10}$ and $\Delta^{15}$, respectively
(see also \cite{Mu} Appendix III Exercise 1.3).
%
\begin{enumerate}
\item[$(1)$] ``Symplectic case".\ $C_n/(\Delta^{10})$ is a finite group.
\item[$(2)$] ``Unitary case".\ $C_n/(\Delta^{15})$ is a finite group.
\end{enumerate}
\noindent
We should remark that $\Delta^{30}=1$ in $C_5$, and $C_5/\Delta^5$ is
a simple group $PSp(4,3)$ (projective symplectic group).
One could try to incorporate deformation of relations
$\Delta^{10}=1$ or $\Delta^{15}=1$ into skein module relations,
but we would speculate that another approach may give appropriate 
skein modules. We will work with oriented links and the following 
useful definition and conjecture \cite{Pr-1,Pr-3,Ki}.
%
\begin{Def}\label{10.3}\begin{enumerate}
\item[$($i$)$]   The $t_n$-move is a local change of an oriented link 
                 which adds $n$ positive half-twists to $L$ $(L_0 \to 
                 L_n)$ as in Figure $16$ $($i$)$.
\item[$($ii$)$]  The $\bar t_k$-move is a local change of an oriented 
                 link which adds $k$ right handed half-twists to $L$ 
                 $(L_0 \to L_k)$ as in Figure $16$ $($ii$)$, where 
                 $k$ is an even integer $($notice anti-parallel 
                 orientation of the strings involved in the move$)$.           
\item[$($iii$)$] Two links are $t_n,\bar t_k$ equivalent if one 
                 can start from the first one and reach the second 
                 one by using $t_n$-moves, $\bar t_k$-moves, and 
                 their inverses.
\end{enumerate}\end{Def}
%
\par\vspace{3mm}\begin{center}\framebox[2cm]{Figure $16$}
\par\vspace{3mm}\end{center}
%
\begin{Conj}\label{10.4}\begin{enumerate}
\item[(i)]  Every oriented link in $S^3$ is $t_3,\bar t_6$ equivalent
            to a trivial link.
\item[(ii)] Every oriented link in $S^3$ is $t_3,\bar t_4$ equivalent
            to a trivial link.
\end{enumerate}\end{Conj}
\noindent
We speculate that there is the following connection 
between our moves and Assion ideals.
%
\begin{problem}\label{10.5}\begin{enumerate}
\item[(i)]  The 5-tangle associated with the 5-braid $\Delta^{10}$ 
            is $3$-move equivalent to the trivial 5-braid tangle.
\item[(ii)] The 5-tangle associated with the 5-braid $\Delta^{15}$ is 
            $t_3,\bar t_4$ equivalent to the trivial 5-braid tangle.
\end{enumerate}\end{problem}
\noindent
Note that from (i) (and \cite{Ch}) it follows immediately that closed 
5-braids are $3$-equivalent to trivial links. B.Wajnryb proved (see 
\cite{Wa-1}, p.694) that every link can be reduced to a trivial link 
by $t_3$ moves and $\Delta^{10}$-moves (Wajnryb's move $\xi$ is 
$t_3$-equivalent to $\Delta^{10}$-move). Thus the positive answer to (i) 
yields the positive answer to the Montesinos-Nakanishi conjecture.
Generally the positive solution to Problem \ref{10.5} would allow 
partial solution to Conjectures \ref{GeneralConjMN} and 
\ref{ConjGenerator}. One should comment, at least shortly, on the 
background of Problem \ref{10.5}. We follow approach and notation of 
\cite{Pr-5} (compare \cite{C-F}, \cite{F-R}). Let $D$ be an oriented 
link diagram. We associate $D$ with the Alexander-Burau module over 
the ring $R$ (with invertible element $t$), as follows. To every arc 
of the diagram we associate a generator (variable) $y_i$ and every 
crossing, $p$, of the diagram yields  the relation 
$(1-t^{\epsilon})y_i +t^{\epsilon}y_k - y_j =0$ where $\epsilon = 
\pm 1$ is the sign of the crossing $p$, see Figure $17$. 
If we think of a neighborhood of a crossing as a $2$-tangle (read 
from left to right), our relation yields a $2 \times 2$ matrix of a 
linear map (we can say that we consider a contra-variant functor from 
the category of 2-tangles to $R$-modules category). For tangles of 
Figure $17$ matrices are: 
%
\[ A_1            = \left[ \begin{array}{cc}  1-t^{-1} & t^{-1} \\
                                                     1 &  0
                           \end{array} \right], \  \ 
   A_2 = A_1^{-1} = \left[ \begin{array}{cc}         0 & 1 \\
                                                     t & 1-t 
                           \end{array} \right]
\] \ \\
%
\par\vspace{3mm}\begin{center}\framebox[2cm]{Figure $17$}
\par\vspace{3mm}\end{center}
%
In particular the matrix corresponding to the move of Figure $18$
is: \par
%
$A_1^3= 
  \left[ \begin{array}{cc}  1-t^{-1} & t^{-1} \\ 1 &  0
           \end{array} \right]^3 =
  \left[ \begin{array}{cc} 1-t^{-3}(t^2-t+1) & t^{-3}(t^2-t+1) \\
                             t^{-2}(t^2-t+1) & 1- t^{-2}(t^2-t+1) 
           \end{array} \right] $ \par\vspace{3mm}\hspace{4mm}
$\equiv \left[ \begin{array}{cc}  1 & 0 \\  0 & 1 
                  \end{array} \right]\ mod (t^2-t+1) $ 
%
\par\vspace{3mm}\begin{center}\framebox[2cm]{Figure $18$}
\par\vspace{3mm}\end{center}
%
Generally we get that the $t_3$-move preserves the Alexander-Burau 
module iff $t^2-t+1=0$ in the ring $R$. We can repeat our analysis 
for $\bar t_{2n}$-moves. In particular the Figure $19$
illustrates an example of $\bar t_2$-moves, with the matrix 
%
\[ B_1 = \left[ \begin{array}{cc}      1-t & t \\
                                         1 & 0 
                  \end{array} \right]
         \left[ \begin{array}{cc} 1-t^{-1} & t^{-1} \\
                                         1 &  0
                  \end{array} \right] =
         \left[ \begin{array}{cc} 2-t^{-1} & t^{-1}-1 \\
                                  1-t^{-1} &  t^{-1}
                  \end{array} \right]
\] 
%
\par\vspace{3mm}\begin{center}\framebox[2cm]{Figure $19$}
\par\vspace{3mm}\end{center}
%
We can now analyze the equivalences $B_1^3 \equiv Id$ and
$B_1^2\equiv Id$ noting that \par
%
\[ B_1^k = \left[ \begin{array}{cc} 1+k(1-t^{-1}) & k(t^{-1}-1) \\
                                      k(1-t^{-1}) & 1-k(1-t^{-1})
\end{array} \right] \]
%
\par\vspace{3mm}\noindent Our results are summarized in the 
following lemma (compare \cite{Pr-1}). 
%
\begin{Lem}\label{}\begin{enumerate}
\item [$(1)$] The Alexander-Burau module $($of an oriented link 
              or a tangle$)$ is unchanged by a $t_3$-move if and 
              only if $t^2-t+1= 0$ in the ring $R$.
\item [$(2)$] The Alexander-Burau module is unchanged by $t_3$ 
              and $\bar t_6$-moves iff $t^2-t+1= 0$ and $3=0$ in $R$.
              Notice that $t^2-t+1= (t+1)^2$ follows from $3=0$.
\item [$(3)$] The Alexander-Burau module is unchanged by $t_3$ 
              and $\bar t_4$-moves iff $t^2-t+1= 0$ and $2=0$ in $R$.
\end{enumerate}\end{Lem}
%
\begin{Cor}\label{10.8}
Let $M_L^{(k)}$ denote the $k$-fold cyclic branched cover of $S^3$ 
$($resp. $D^3)$ with the branching set the link $($resp. tangle$)$ 
$L$. Then:   \begin{enumerate}
\item [$(1)$]  $H_1(M_L^{(2)},Z_3)$ is unchanged by a $t_3$ 
               and $\bar t_6$-moves $($more generally by 3-moves$)$.
\item [$(2)$]  $H_1(M_L^{(3)},Z_2)$ is unchanged by a $t_3$ and 
               $\bar t_4$-moves.
\end{enumerate}\end{Cor}
%
For our heuristic argument for a positive answer to Conjecture 
\ref{10.4} it remains to show that $\Delta^{10}$-move preserves 
Alexander-Burau matrix modulo $(t+1,3)$ and $\Delta^{15}$-move 
preserves Alexander-Burau matrix modulo $(t^2-t+1,2)$. It is a 
tedious but easy task. In particular Figure $20$ illustrates the fact 
that the matrix $M(\Delta^{-1})$ for $\Delta^{-1}= 
\sigma_4^{-1}\sigma_3^{-1}\sigma_2^{-1}\sigma_1^{-1}$ is
%
\[ M(\Delta^{-1}) = \left[ \begin{array}{ccccc}
 0 & 0 & 0 & 0 & 1 \\
 t & 0 & 0 & 0 & 1-t \\
 0 & t & 0 & 0 & 1-t \\
 0 & 0 & t & 0 & 1-t \\ 
 0 & 0 & 0 & t & 1-t 
\end{array} \right]
\] 
%
\par\vspace{3mm}\begin{center}\framebox[2cm]{Figure $20$}
\par\vspace{3mm}\end{center}
%
We get modulo $(t^2-t+1)$:
\[ M(\Delta^{-10})\equiv \left[ \begin{array}{ccccc}
 -2& -1 + 2t & 1 + t & 2 - t & 1 - 2t\\
 -2 + t & -1 + t & 1 + t & 2 - t & 1 - 2t\\
 -2 + t & -1 + 2t & 1 & 2 - t & 1 - 2t\\
 -2 + t & -1 + 2t & 1 + t &  2 - 2t & 1 - 2t\\
 -2 + t & -1 + 2t & 1 + t & 2 - t & 1 - 3t 
\end{array} \right] \]
%
which is the identity matrix iff $3\equiv 0$ and $t+1 \equiv 0$.
In fact we have more generally: $M(\Delta^{-10}) \equiv \ Id \ mod(t+1)$.
One should stress here that the $\Delta^{-10}$ move is not
preserving the Alexander-Burau module modulo $(t^2-t+1,3)$
therefore this move is not a combination of $t_3$ and $\bar t_6$
moves (we conjecture only that it is a combination of 3-moves).
We have also (still modulo $(t^2-t+1)$):
%
\[ M(\Delta^{-15}) \equiv 
\left[ \begin{array}{ccccc} 
 -3 + 2t & 2t & 2 & 2 - 2t & -2t\\
 -2 + 2t & -1 + 2t & 2 & 2 - 2t & -2t\\
 -2 + 2t & 2t & 1 & 2 - 2t & -2t\\
 -2 + 2t & 2t & 2 & 1 - 2t & -2t\\
 -2 + 2t & 2t & 2 & 2 - 2t & -1 - 2t
\end{array} \right]
\]
%
which is the identity matrix iff $2\equiv 0$.
It follows from Coxeter work that $M(\Delta^{-30}) \equiv 
\ Id \ mod \ (t^2-t+1)$. We have more generally: 
$M(\Delta^{-30})\equiv \ Id \ mod (t^3+1)$.
%
\section*{Acknowledgment}
We would like to thank Bruce Westbury for very helpful discussion.
%
%

%
\addtocounter{figure}{-20}
\par\vspace{3mm}\begin{figure}[htbp]\begin{center}
\includegraphics[trim=0mm 0mm 0mm 0mm, width=.7\linewidth]
{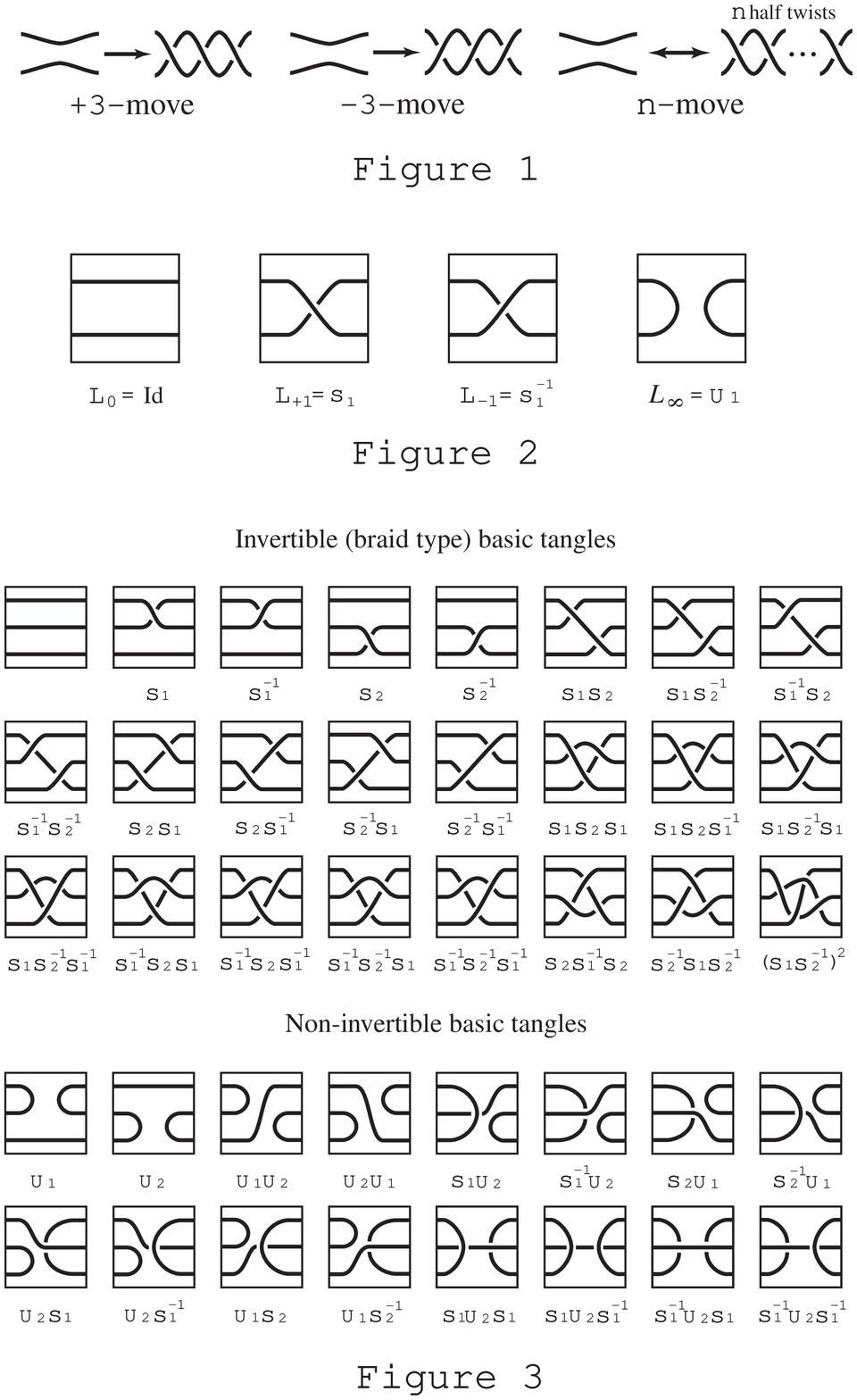}\end{center}\end{figure}\par\vspace{2mm}
%
\par\vspace{3mm}\begin{figure}[htbp]\begin{center}
\includegraphics[trim=0mm 0mm 0mm 0mm, width=.7\linewidth]
{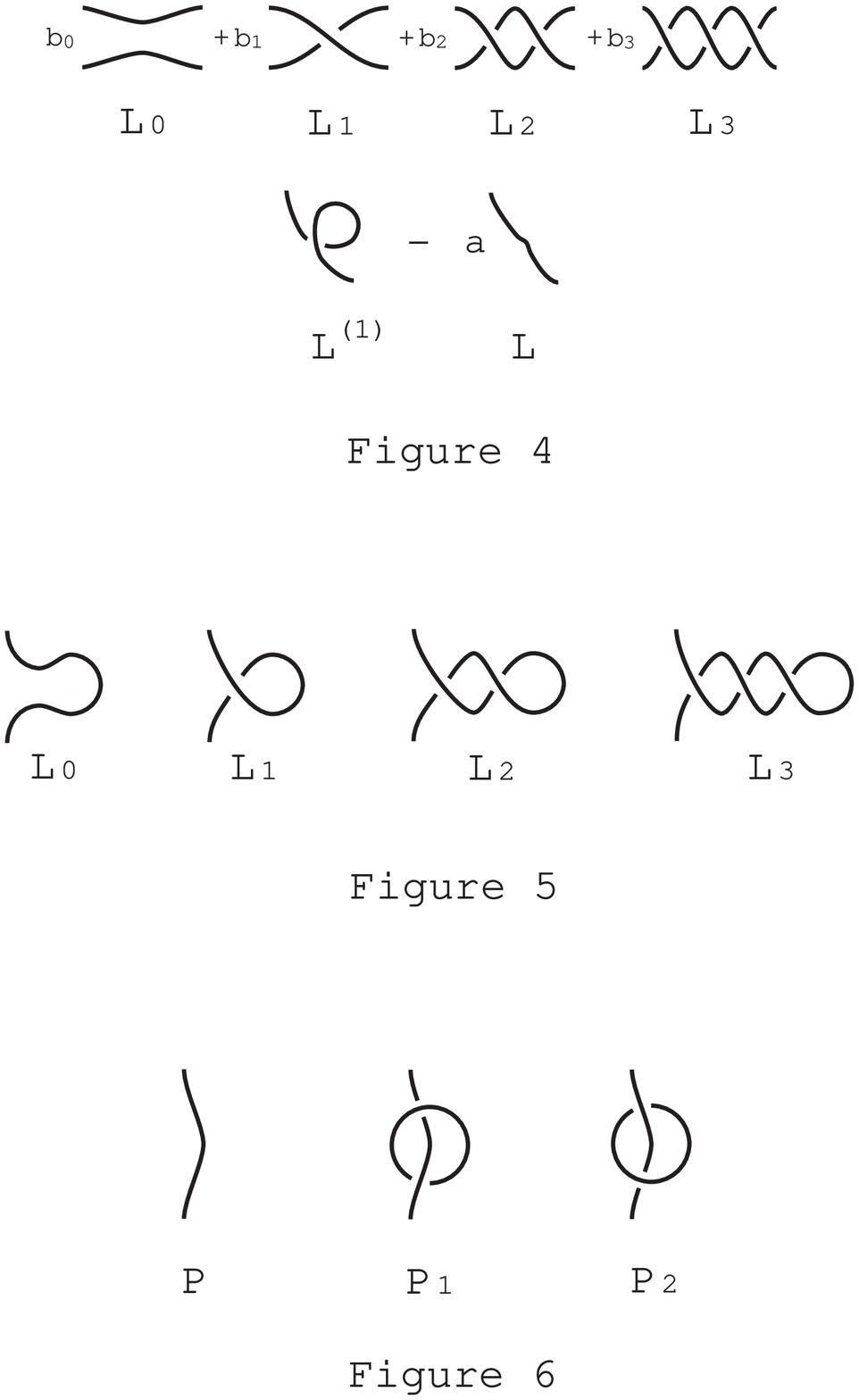}\end{center}\end{figure}\par\vspace{2mm}
%
\par\vspace{3mm}\begin{figure}[htbp]\begin{center}
\includegraphics[trim=0mm 0mm 0mm 0mm, width=.7\linewidth]
{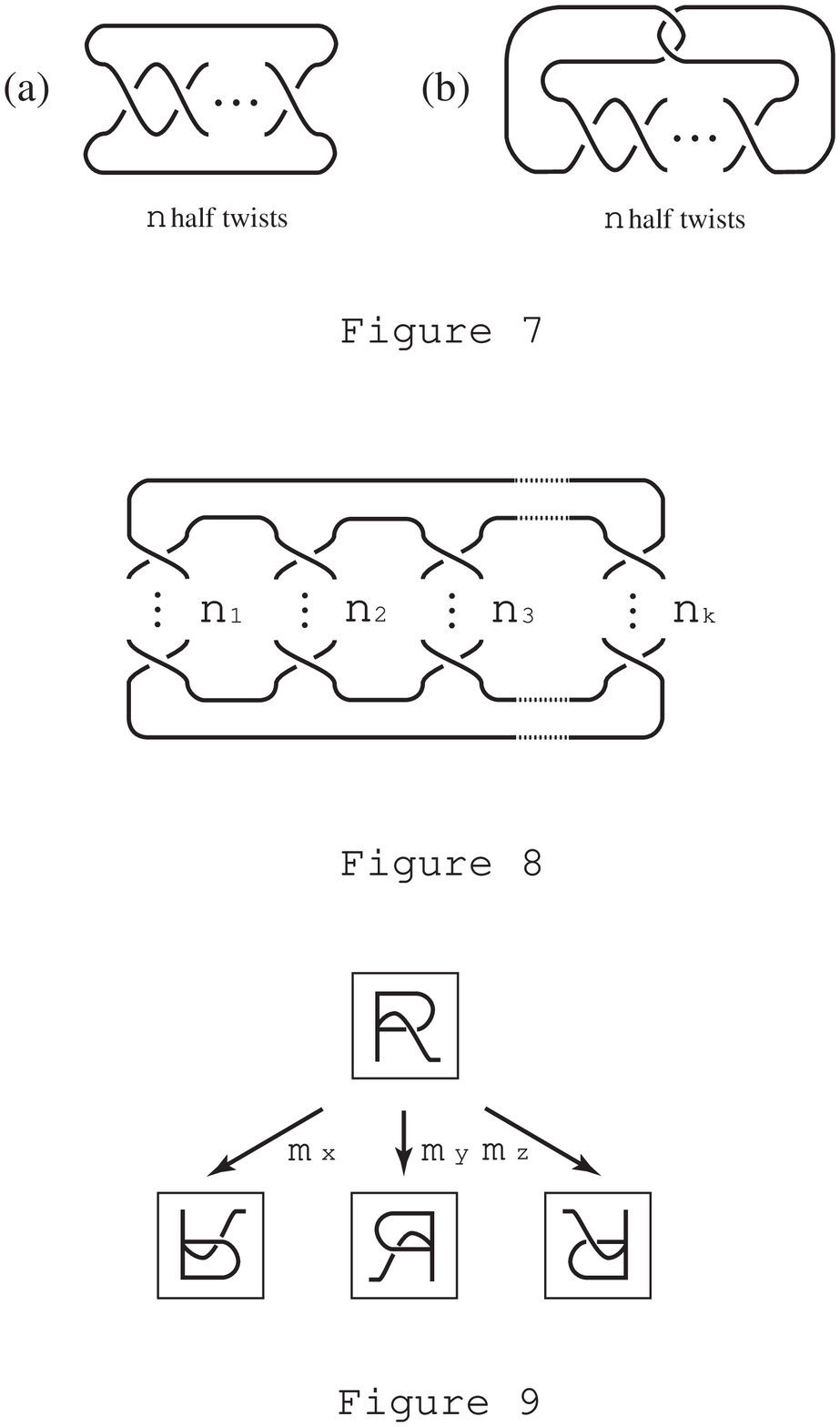}\end{center}\end{figure}\par\vspace{2mm}
%
\par\vspace{3mm}\begin{figure}[htbp]\begin{center}
\includegraphics[trim=0mm 0mm 0mm 0mm, width=.6\linewidth]
{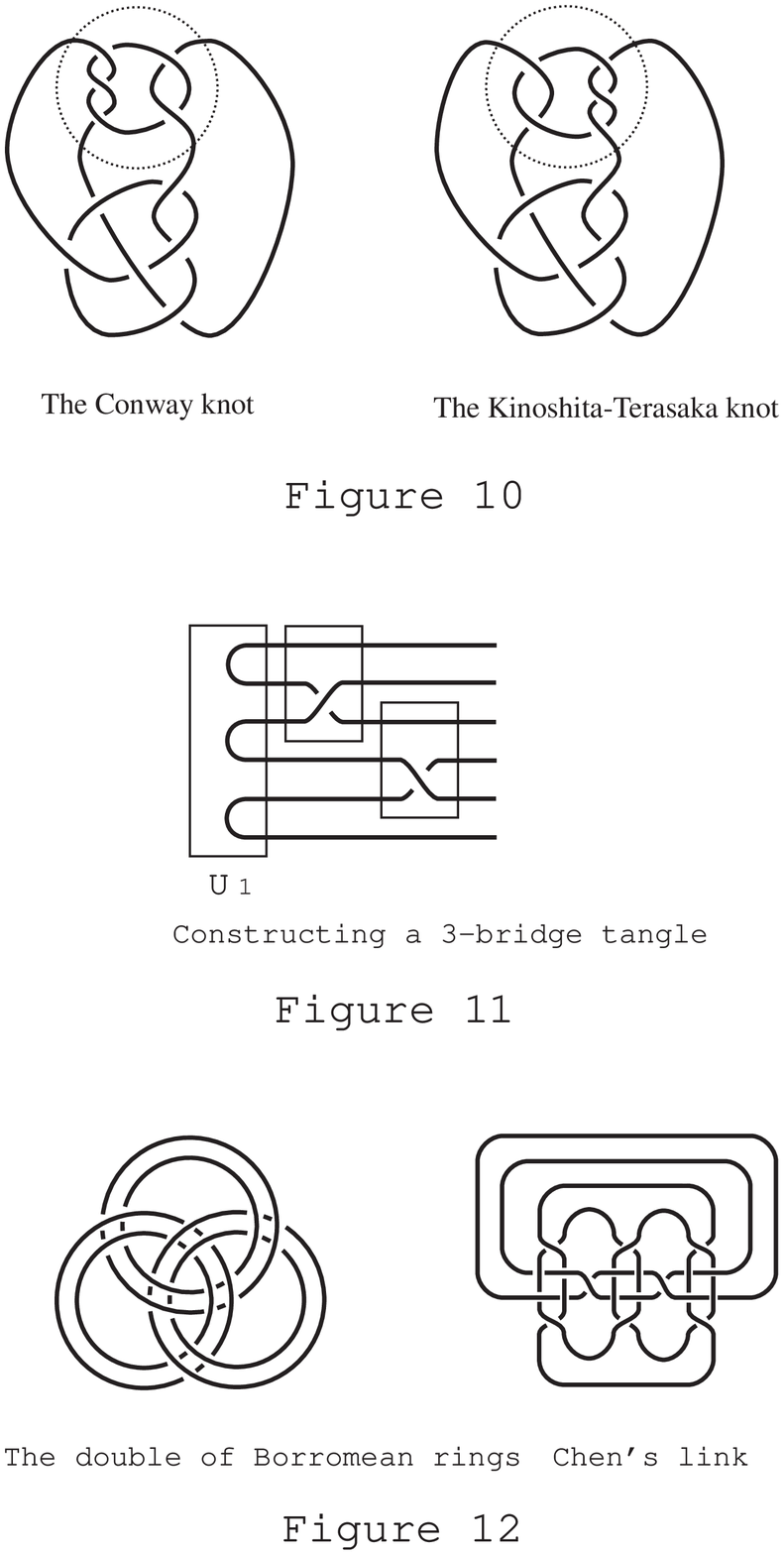}\end{center}\end{figure}\par\vspace{2mm}
%
\par\vspace{3mm}\begin{figure}[htbp]\begin{center}
\includegraphics[trim=0mm 0mm 0mm 0mm, width=.8\linewidth]
{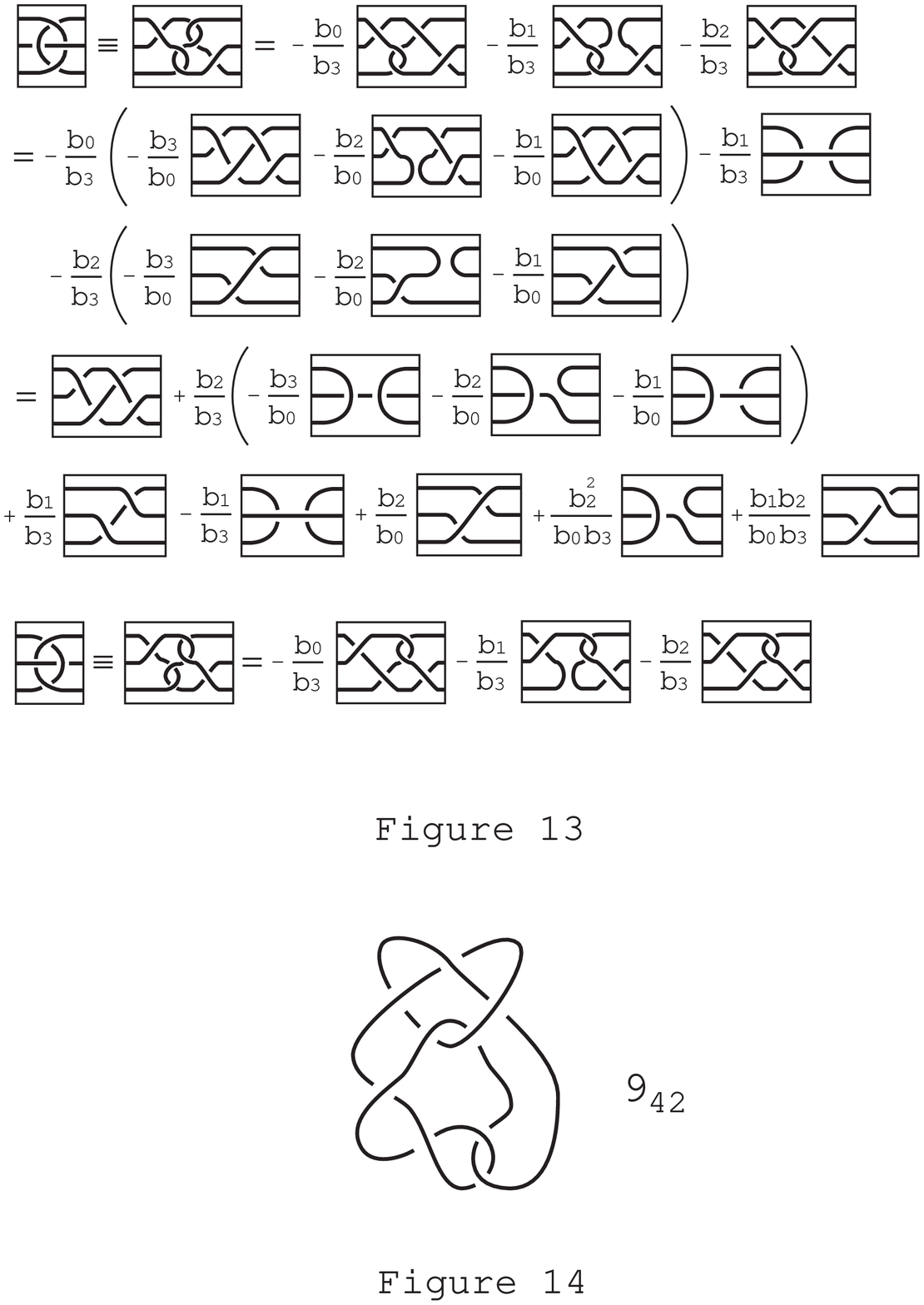}\end{center}\end{figure}\par\vspace{2mm}
%
\par\vspace{3mm}\begin{figure}[htbp]\begin{center}
\includegraphics[trim=0mm 0mm 0mm 0mm, width=.8\linewidth]
{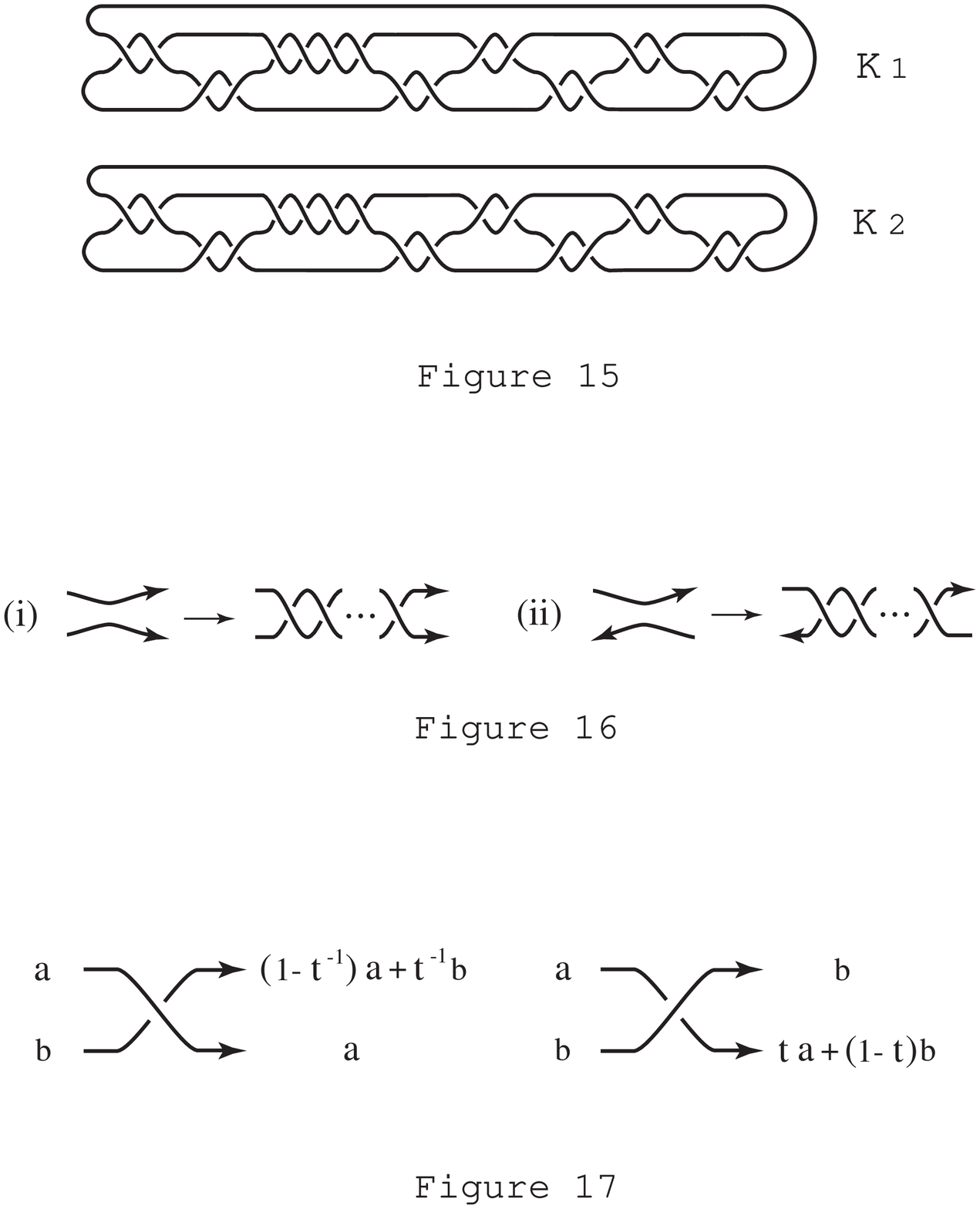}\end{center}\end{figure}\par\vspace{2mm}
%
\par\vspace{3mm}\begin{figure}[htbp]\begin{center}
\includegraphics[trim=0mm 0mm 0mm 0mm, width=.8\linewidth]
{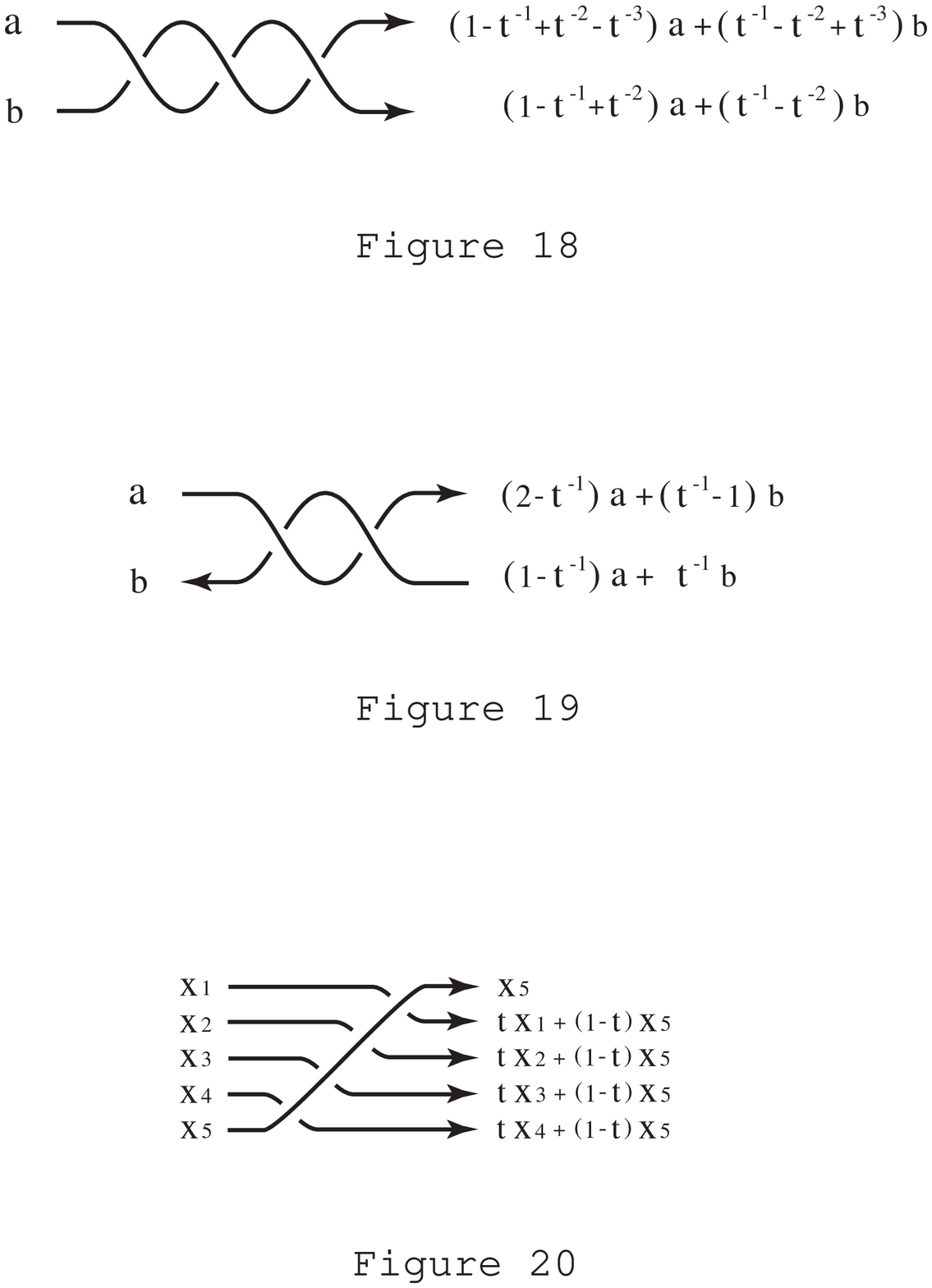}\end{center}\end{figure}\par\vspace{2mm}
%
\end{document}